# ARCHIMEDE, THEON de SMYRNE et √3

Abd Raouf Chouikha, Université Paris-Sorbonne, Paris-Nord

**Résumé :** *On sait que l'algorithme de Théon de Smyrne (70-135 ap JC) a permis de mettre en évidence des encadrements fins de √2 par des rationnels. Pourtant ce même algorithme s'applique aussi à √3 et permet de retrouver les fameuses inégalités d'Archimède. Une question intéressante est de savoir si cette méthode itérative très simple exposée par Théon n'était pas connue des grecs avant lui, notamment par Archimède ou ses contemporains. Sachant que Théon de Smyrne avait déjà compilé des travaux antérieurs. Notons que cet algorithme avait permis aussi d'ouvrir une brèche vers les équation de Pell-Fermat et les fractions continues.*

**Abstract :** *We know that the algorithm of Theon of Smyrna (70-135 AD) made it possible to highlight fine frames of √2 by rationals. However, this same algorithm also applies to √3 and makes it possible to find the famous Archimedes inequalities. An interesting question is whether this very simple iterative method exposed by Theon was not known to the Greeks before him, notably by Archimedes or his contemporaries. Knowing that Theon of Smyrna had already compiled previous works. Note that this algorithm had also made it possible to open a breach towards the Pell-Fermat equations and continued fractions.*

## I – Origine et influence de l'algorithme de Théon de Smyrne

### 1. Description de l'algorithme de Théon de Smyrne [6]

Considérons une suite de carrés emboîtés de côté $x_0, x_1, x_2, ... x_n$ obtenus en ajoutant à la mesure du côté précédent la mesure de sa diagonale. Les nombres $x_n$ sont dits *latéraux*. Soit $y_n$ la diagonale du carré de côté $x_n$. Les $y_n$ sont dits *diagonaux*.

On a $x_{n+1} = x_n + y_n$ et selon le théorème de Pythagore :

$$y_n = \sqrt{2}\ x_n \quad \text{et} \quad y_{n+1} = \sqrt{2}\ x_{n+1} = \sqrt{2}\ (x_n + y_n)$$

Par conséquent : $y_{n+1} = y_n + 2 x_n \quad \text{et} \quad x_{n+1} = x_n + y_n$

Posons alors $r_n = y_n/x_n$. Il est aisé de montrer :

$$r_{n+1} = \frac{r_n + 2}{r_n + 1}$$

On construit ainsi une <u>suite récurrente</u> de la forme $r_{n+1} = f(r_n)$. Ce qui signifie si $r_n$ est une approximation de $\sqrt{2}$, alors $r_{n+1}$ est une meilleure approximation. Il est facile de prouver la convergence vers $\sqrt{2}$ en choisissant par exemple $r_0 = 1$, ce qui correspond à choisir $x_0 = 1$ et $y_0 = 1$. On a donc

$$r_{n+1} - \sqrt{2} = (\sqrt{2} - 1)\frac{(\sqrt{2} - r_n)}{(r_n + 1)} \quad \Rightarrow$$

$$|r_{n+1} - \sqrt{2}| < \left(\tfrac{1}{2}\right) |\sqrt{2} - r_n| \quad \Rightarrow \quad |r_n - \sqrt{2}| < \left(\tfrac{1}{2}\right)^n [\sqrt{2} - 1].$$

## 2. Origine de l'Algorithme de Théon de Smyrne, [7]

- Théon de Smyrne (I[er]–II[e] siècle apr. J.-C.) est un mathématicien et philosophe grec de l'Antiquité tardive. Son algorithme pour $\sqrt{2}$ est décrit dans son ouvrage *Expositio rerum mathematicarum* (ou Des connaissances mathématiques utiles à la lecture de Platon), où il explique des méthodes arithmétiques et géométriques. «*Je ferai un résumé et une esquisse concise des théorèmes mathématiques qui sont particulièrement nécessaires aux lecteurs de Platon.* » dit-il dans l'introduction.
- A aucun moment, Théon ne revendique des résultats ou prétend inventer de nouveaux théorèmes. Il expose simplement les mathématiques que l'on doit connaître afin de lire Platon.
- Son but : Fournir des approximations de $\sqrt{2}$ pour des applications pratiques (géométrie, architecture, astronomie) et illustrer les propriétés des nombres irrationnels, déjà connus depuis les Pythagoriciens (V[e] siècle av. J.-C.).

**Méthode d'antiphérèse (soustraction réciproque)**

Le néo-pythogoricien Théon de Smyrne (70~135), inspiré de la méthode d'antiphérèse d'Euclide, a utilisé le principe de suites encadrantes pour approcher la valeur de $\sqrt{2}$ (a/b, (a+2b)/(a+b) ... : 3/2, 7/5, 17/12 ...). Diophante d'Alexandrie (214~298) aura <u>besoin</u> plus tard d'une <u>théorie des nombres</u> unifiant les entiers, les rationnels et les irrationnels, pour la résolution d'une <u>équation</u> diophantienne.

# 3. Texte de Théon de Smyrne sur les nombres latéraux et diagonaux, [4]

A titre d'information nous donnons la traduction du texte original de Théon. On pourra apprécier la qualité de la pédagogie, la précision de la description de l'algorithme, un peu pénible parfois à lire à cause du manque de notations symboliques :

[Exposition des connaissances mathématiques utiles pour la lecture de Platon / Théon de Smyrne,... ; traduite pour la première fois du grec en français par J. Dupuis](#)

### Des nombres latéraux et des nombres diagonaux

XXXI. De même que les nombres ont en puissance les rapports des triangulaires, des tétragones, des pentagones et des autres figures, de même nous trouverons que les rapports des nombres latéraux et des nombres diagonaux se manifestent dans les nombres selon des raisons génératrices, car ce sont les nombres qui harmonisent les figures. Donc comme l'unité est le principe de toutes les figures, selon la raison suprême et génératrice, de même aussi le rapport de la diagonale et du côté se trouve dans l'unité.

Supposons par exemple deux unités dont l'une soit la diagonale et l'autre le côté, car il faut que l'unité qui est le principe de tout soit en puissance le côté et la diagonale; ajoutons au côté la diagonale et à la diagonale ajoutons deux côtés, car ce que le côté peut deux fois, la diagonale le peut une fois. Dès lors la diagonale est devenue plus grande et le côté plus petit. Or, pour le premier côté et la première diagonale, le carré de la diagonale unité sera moindre d'une unité que le double carré du côté unité, car les unités sont en égalité, mais un est moindre d'une unité que le double de l'unité. Ajoutons maintenant la diagonale au côté, c'est- à-dire une unité à l'unité, le côté vaudra alors 2 unités; mais, si nous ajoutons deux côtés à la diagonale, c'est-à-dire 2 unités à l'unité, la diagonale vaudra 3 unités; le carré construit sur le côté 2 et 4, et le carré de la diagonale est 9 qui est plus grand d'une unité que le double carré de 2.

De même ajoutons au côté 2 la diagonale 3, le côté deviendra 5. Si à la diagonale 3 nous ajoutons deux côtés, c'est-à-dire 2 fois 2, nous aurons 7 unités. Le carré construit sur le côté est 2, et celui qui est construit sur la diagonale 7 est 49, qui est moindre d'une unité que le double 50 du carré 25. De nouveau, si au côté 5 on ajoute la diagonale 7, on obtient 12 unités; et si à la diagonale 7 on ajoute 2 fois le côté 5, on aura 17 dont le carré 289) est plus grand d'une unité que le double (288) du carré de 12. Et ainsi de suite en continuant l'addition. La proportion alterne : le carré construit sur la diagonale sera tantôt plus petit, tantôt plus grand, d'une unité, que le double carré construit sur le côté, en sorte que ces diagonales et ces côtés seront toujours exprimables.

Inversement les diagonales comparées aux côtés, en puissance, sont tantôt plus grandes d'une unité que les doubles, tantôt plus petites d'une unité. Toutes les diagonales sont donc, par rapport aux carrés des côtés, doubles alternativement par excès et par défaut, la même unité combinée également avec tous, rétablissant l'égalité, en sorte que le double ne pèche ni par excès, ni par défaut; en effet, ce qui manque dans la diagonale précédente se trouve en excès, en puissance, dans la diagonale qui suit.

Cela veut dire en notation moderne, qu'on définit une suite de nombres rationnels :

$r_0 = 1$, $r_1 = 3/2$, $r_2 = 7/5$, $r_3 = 17/12$, $r_4 = 41/29$, $r_5 = 99/70$, $r_5 = 239/169$, ….

On construit donc une suite d'encadrements emboîtés de $\sqrt{2}$ de plus en plus fins. On établit ainsi $(r_{2n})$ et $(r_{2n+1})$ 2 suites adjacentes convergeant vers $\sqrt{2}$.

## 4. Proclus et les equations de Pell, [9,10]

Proclus, philosophe et mathématicien grec du Ve siècle, est connu pour ses observations sur la relation entre les approximations de $\sqrt{2}$. Son travail s'appuie sur des études pythagoriciennes antérieures sur les «côtés et diamètres », qui ont généré des approximations rationnelles pour $\sqrt{2}$ à travers des suites telles que 1/1, 3/2, 7/5, 17/12, etc.

Proclus a observé que ces approximations obéissaient à l'équation $x^2 - 2y^2 = \pm 1$. Par exemple, la paire $x = 7$, $y = 5$ satisfait $7^2 - 2 \times 5^2 = -1$, et $x = 17$,

y =12 donne 17²−2×12² =1,…. Cette idée a mis en évidence comment les solutions itératives pouvaient approximer des nombres irrationnels comme √2, une idée fondamentale formalisée plus tard dans l'étude des équations de Pell.

Bien que Proclus n'ait pas développé de méthode générale pour résoudre les équations de Pell, sa reconnaissance de ce modèle a contribué à une meilleure compréhension des équations diophantiennes quadratiques dans l'Antiquité. Son travail illustre les premiers efforts grecs pour relier la théorie des nombres aux problèmes géométriques, en particulier ceux impliquant des quantités irrationnelles. L'étude systématique des équations de Pell allait cependant se développer plus tard dans les mathématiques indiennes avec la méthode chakravala de Brahmagupta et en Europe au cours du XVIIe siècle.

Proclus pensait que les Pythagoriciens calculaient **√2** en utilisant le fait suivant, (traduit en notation moderne) : si $a_n$, $b_n$ satisfont $(a_n)^2 - 2(b_n)^2 = \pm 1$ et posons $a_{n+1} = a_n + 2 b_n$, $b_{n+1} = a_n + b_n$ alors

$$(a_{n+1})^2 - 2(b_{n+1})^2 = \mp 1.$$

En partant de $a_1 = b_1 = 1$ cela génère une séquence de bornes supérieures et inférieures qui s'améliorent progressivement en alternance. Il s'agit en fait d'une procédure systématique pour résoudre l'équation de Pell lorsque c = 2 et m = ±1. Cela rappelle curieusement l'algorithme de Théon.

$$1 < 7/5 < 41/29 < 239/169 \ldots < \sqrt{2} < \ldots < 99/70 < 17/12 < 3/2 < 2$$

**Rappelons l'équation de Pell-Fermat pour n entier, [10]**

$$y^2 - n x^2 = \pm 1$$

Cette équation a une infinité de <u>solutions</u> si *n* n'est pas un carré (n'a pas de facteur carré).

- Dénommée d'abord par <u>Euler</u> (1707-1783), du nom du mathématicien anglais <u>Pell</u> (1610-1685).
- <u>Bhaskara</u> (600-680) a tenté de résoudre cette équation.
- Le mathématicien indien <u>Brahmagupta</u> (598-670) est le premier à avoir décrit l'ensemble des solutions de cette équation.

- Fermat (1601-1665) conjecture qu'elle a une infinité de solutions.
- Lagrange (1736-1813) le prouve un siècle plus tard.
- Brouncker (1620-1684) trouve les relations qui relient la fraction continue d'un nombre quadratique à l'équation de Pell-Fermat.

Fraction continue de racine de 2:

$$\sqrt{2} = [1; 2, 2, 2 \ldots]$$

## 5. Cet algorithme était-il connu avant Théon ? [5,13]

- Les Pythagoriciens connaissaient l'irrationalité de $\sqrt{2}$ dès le V$^e$ siècle av. J.-C., mais ils n'ont pas laissé de traces écrites explicites de méthodes d'approximation (leur savoir était souvent secret).
- Archimède (III$^e$ siècle av. J.-C.) et d'autres mathématiciens grecs utilisaient des approximations fractionnaires pour $\sqrt{2}$, comme $99/70 \approx 1{,}4142857$ (proche de $\sqrt{2} \approx 1{,}4142135$). Ces approximations très simples en l'occurrence ne suggèrent-elles pas une connaissance empirique de méthodes itératives ? Au point qu'Archimède aurait pu l'appliquer en d'autres situations.
- La méthode de Théon elle-même pourrait s'inspirer de techniques plus anciennes, mais c'est bien lui qui l'a formalisée sous forme explicite d'un algorithme itératif. Cependant, n'y a pas de preuve directe que cette méthode précise était utilisée avant lui.
- Théon a-t-il inventé cette méthode ? Il est difficile de l'affirmer, car les sources antérieures sont fragmentaires. Nous savons que ses ouvrages sont d'une qualité pédagogique indéniable. Cependant, il ne revendique pas ces résultats, mais il est le premier à l'avoir clairement décrite et systématisée dans un texte.
- Toutefois **P. Tannery**, [13,p.150] disait à ce sujet : « …Mais il faut remarquer que Théon est lui-même un compilateur, ce dont il ne se cache nullement…. »
- Les anciens Grecs connaissaient des approximations de $\sqrt{2}$ et utilisaient des raisonnements géométriques (comme la diagonale du carré), mais la méthode itérative proposée par Théon représente une avancée dans la formalisation algébrique.

### a - La diagonale du carré

Les Grecs ont également étudié √2 dans un contexte géométrique. Selon le **théorème de Pythagore**, la diagonale *d* d'un carré de côté 1 est donnée par :

$$\boldsymbol{d = \sqrt{(1^2+1^2)} = \sqrt{2}}$$

Cela a conduit à des constructions géométriques pour approximer √2, comme l'utilisation de la spirale de Théodore ou d'autres méthodes basées sur des proportions.

### b - Tablettes babyloniennes, [15]

Bien que les Grecs aient développé leurs propres méthodes, ils ont probablement été influencés par les Babyloniens, qui avaient déjà des approximations précises de √2. Par exemple, une tablette babylonienne (YBC 7289) donne une approximation de √2 en base 60 :

$$\sqrt{2} \approx 1 + 24/60 + 51/60^2 + 10/60^3 \approx 1{,}41421296$$

Cette valeur est remarquablement proche de la valeur réelle (1,41421356...).

### c - Connexions avec les triplets pythagoriciens, [14]

On a la relation pour ***n* impair** : $x_n^2 - 2 y_n^2 = -1$

$x_{n+1}^2 - 2 y_{n+1}^2 = (x_n + 2 y_n)^2 - 2(x_n + y_n)^2 = -(x_n^2 - 2 y_n^2) = -1 = (-1)^{n+1}$

donne un triplet pythagoricien : $((x_n - 1)/2)^2 + ((x_n + 1)/2)^2 = y_n^2$

(remarquons que $x_n$ est nécessairement impair)

Pour ***n* pair**, $x_n^2 - 2 y_n^2 = 1$ : $(x_n^2 - 1)/2 = y_n^2$ est un carré parfait.

Rappelons que Proclus attribue à Pythagore la découverte de la formule générale des triplets où $n > 0$ [Triplet pythagoricien — Wikipédia](#) :

$$(2n + 1,\ 2n^2 + 2n,\ 2n^2 + 2n + 1)$$

# II - Généralisation de l'algorithme de Théon

On peut généraliser l'algorithme à la racine carrée d'un nombre positif quelconque :

Soit A un nombre positif; on pose $x_0 = 1$, $y_0 = 1$, et, pour tout entier naturel n :

$$x_{n+1} = x_n + y_n \quad \text{et} \quad y_{n+1} = y_n + A x_n$$

Un calcul donne après avoir posé

$r_n = y_n / x_n$

$$r_{n+1} = \frac{r_n + A}{r_n + 1}$$

Si A est un carré parfait $A = a^2$, a entier naturel alors on a une suite constante de terme tous égaux à a :

$$r_n = a, \text{ pour tout } n$$

## 1. Cas où A=3

Contexte géométrique et théorique

Les Grecs associaient √3 à des figures géométriques, comme le **triangle équilatéral** (où la hauteur se déduit par $h = \sqrt{3}/2 \times$ côté´). Cette approche pratique, combinée à leur compréhension des nombres irrationnels, les incitait à développer des approximations utiles pour l'architecture ou pour les calculs astronomiques.

On ne connaît pas la méthode qu'il a utilisée pour arriver à ce résultat. Davies [2] spécule sur plusieurs possibilités dans son article, et conclut qu'Archimède a probablement utilisé la méthode de Héron, basée sur un algorithme de calcul de racines carrées qui repose sur le fait que la moyenne arithmétique de deux nombres est supérieure à leur moyenne géométrique, qui est elle-même supérieure à leur moyenne harmonique, (Davies, [2] 2011).

## 2. Algorithme de Héron d'Alexandrie (I[er] siècle), [6]

Héron a généralisé une méthode itérative pour calculer les racines carrées, applicable à √3. Son algorithme consiste à :

1. Choisir une estimation initiale $x_0$ (par exemple, 2).
2. Calculer la moyenne entre $x_n$ et $3/x_n$ pour obtenir $x_{n+1}$ :

$$x_{n+1} = (\tfrac{1}{2})(x_n + 3/x_n)$$

Exemple pour $\sqrt{3}$ :

- $x_1 = (2+3/2)/2 = 7/4$

  - Itération 1 : $x_1 = (2+3/2)/2 = 1{,}75$

  $x_2 = 1{,}75 + 3/1{,}752 = 97/56$

  - Itération 2 : $x_2 = (7/4 + 3.4/7)/2 \approx 1{,}732142$, déjà très proche de la valeur

    réelle.

Cette méthode, redécouverte plus tard en Europe, permettait d'atteindre une précision remarquable en quelques étapes.

Cependant, Sondheimer [12] avait remarqué qu'un algorithme explicite de fractions continues aurait été difficile à réaliser pour les Grecs en raison de toutes les longues divisions requises (Sondheimer et Rogerson, Numbers and Infinity, 1981).

## 3. Utilisation de l'algorithme de Théon, [2]

Utilisons maintenant la méthode de Théon. Pour A = 3, on trouve :

$r_0 = 1$, $r_1 = 2$, $r_3 = 5/3$, $r_4 = 7/4$, $r_5 = 19/11$, $r_6 = 26/15$, $r_7 = 71/41$, $r_8 = 97/56$,

**$r_9 = 265/153$**, $r_{10} = 362/209$, $r_{11} = 989/571$, **$r_{12} = 1351/780$**,….

On a donc les encadrements suivants :

$$1 < 5/3 < 19/11 < 71/41 < \mathbf{265/153} < \ldots < \sqrt{3} <$$
$$\ldots \mathbf{1351/780} < 362/209 < 97/56 < 26/15 < 7/4 < 2.$$

On retrouve la fameuse inégalité d'Archimède, utilisée pour donner l'encadrement de Pi :

$$r_9 = \mathbf{265/153} < \sqrt{3} < r_{12} = \mathbf{1351/780} \qquad \Leftrightarrow$$

$$26 - 1/52 < 15\sqrt{3} < 26 - 1/51.$$

On a en plus une estimation de l'erreur

$$r_{n+1} - \sqrt{3} = (\sqrt{3} - 1)(\sqrt{3} - r_n)/(r_n) \Rightarrow$$

$$|r_{n+1} - \sqrt{3}| < (1/2)|\sqrt{3} - r_n| \Rightarrow$$

$$|r_n - \sqrt{3}| < (1/2)^n [\sqrt{3} - 1].$$

On ne connaît pas exactement l'argument utilisé par Archimède pour justifier cette double inégalité mentionnée dans « Mesure du cercle » :

**$265/153 < \sqrt{3} < 1351/780$**

On peut supposer que les méthodes itératives étaient familières à Archimède. Puisqu'il connaissait par exemple 99/70 comme approximation de √2.

D'après Davies [2] : "a good candidate should satisfy the following criteria.

• It should only use methods that were known at the time of Archimedes;

• It should be possible to replace 3 by any other positive integer;

• It should be as short and elementary as possible;

• It should not involve extraordinary ingenuity or tricks."

L'algorithme de Théon de Smyrne ne répond-il pas à tous ces standards ?

Déjà Heath faisait remarquer [6] :

...the calculation [of π] starts from a greater and lesser limit to the value of √3, which Archimedes assumes without remark as known, namely 265/153 < √3 < 1351/780. How did Archimedes arrive at this particular approximation? No puzzle has exercised more fascination upon writers interested in the history of mathematics... Another suggestion...is that the successive solutions in integers of the equations x2-3y2=1 and x2-3y2=-2 may have been found...in a similar way to...the Pythagoreans. The rest of the suggestions amount for the most part to the use of the method of continued fractions more or less disguised.
—T. Heath, *A History of Greek Mathematics*, 1921

**Méthode des moyennes itératives** (Procédé archimédien — Wikipédia, [11] )

Bien qu'Archimède n'ait pas formalisé d'algorithme comme Héron, ses travaux préfigurent des méthodes utilisant des **moyennes harmoniques et géométriques**. Le "procédé archimédien" décrit dans les sources modernes consiste à générer deux suites convergentes vers une limite commune, comme la moyenne arithmético-géométrique. Pour $\sqrt{3}$, cela pourrait correspondre à une initialisation avec des bornes simples (ex. $5/3 \approx 1{,}666$) et des ajustements successifs

Note_on_generalized_Archimedes-Borchardt_algorithm.pdf , [8]

S0025-5718-1984-0725993-3.pdf , [3]

### 4. Fractions continues et rectangles, [5]

Les sources suggèrent qu'Archimède aurait pu utiliser une **représentation géométrique des fractions continues** pour approximer $\sqrt{3}$. En construisant des rectangles dont les rapports de côtés se rapprochent de $\sqrt{3}$, il obtenait des fractions comme $265/153$ ($\approx 1{,}732$) ou $1351/780$ ($\approx 1{,}73205$), bien que ces valeurs spécifiques ne soient pas explicitement citées dans les résultats de recherche

**Notations :**

$$\sqrt{3} = [1, 1, 2, 1, 2, 1, 2, 1, 2...]$$

$$\sqrt{3} = 1 + \cfrac{1}{2 + \cfrac{1}{1 + \cfrac{1}{2 + \cdots}}}$$

### 5. Influences et applications

Archimède a bien influencé des méthodes ultérieures comme l'algorithme de Héron ou l'approche de Borchardt, qui généralisent l'utilisation des moyennes pour les racines carrées. Son encadrement de $\sqrt{3}$, combiné à ses travaux sur $\pi$, illustre sa maîtrise des techniques d'approximation géométrique, posant les bases de l'analyse numérique.

[Procédé archimédien — Wikipédia](), [11]

On sait, grâce aux tablettes d'argile cunéiformes, que les Babyloniens étaient capables de calculer les racines carrées avec précision bien avant l'époque d'Archimède, mais il n'existe aucune preuve contemporaine de la manière dont ils y parvenaient. Au cinquième siècle après J.-C., Proclus affirmait que les Pythagoriciens calculaient √2 en utilisant le fait suivant réécrit en notation moderne. Si $a_n$, $b_n$ satisfont $a_n^2 - 2 b_n^2 = \pm 1$ et que l'on pose : $a_{n+1} = a_n + 2 b_n$, $b_{n+1} = a_n + b_n$ alors $a_{n+1}^2 - 2 b_{n+1}^2 = \mp 1$. En partant de $a_1 = b_1 = 1$, cela génère une séquence de bornes supérieures et inférieures qui s'améliorent progressivement en alternance. Il s'agit en fait d'une procédure systématique pour résoudre l'équation de Pell lorsque c = 2 et m = ±1. Cela rappelle curieusement l'algorithme de Théon.

## **Conclusion** :

Archimède a vraisemblablement approché √3 en utilisant une méthode itérative ou algorithmique assez proche à celle décrite par Théon de Smyrne. Il a aussi combiné intuition géométrique, encadrement par polygones et manipulation de fractions, sans disposer apparemment des outils algébriques modernes. Bien que les détails exacts de sa méthode restent partiellement obscurs, son influence sur le développement des algorithmes numériques est indéniable. Ces travaux illustrent sa capacité à relier abstraction mathématique et applications pratiques, influençant des siècles de recherche en analyse numérique.

D'après **Bernard Beauzamy [1]:** « Beaucoup de gens pensent que la science grecque dans son ensemble est aujourd'hui dépassée : nous nous réjouissons d'être capables de calculer des milliards de décimales de π en quelques millisecondes, et méprisons les Grecs qui, pour nous, n'étaient capables que d'une approximation au moyen de bouts de bois.

Citons E.T. Bell, dans son gros ouvrage "Men of Mathematics" : "It is only possible that Archimedes, could he come to life long enough to take a post-graduate course in mathematics and physics, would understand Einstein, Bohr, Heisenberg and Dirac better than they would understand themselves."

….. S'il est permis de parler de durée de retour pour le génie, on dira que Gauss est de classe 500 (on voit un Gauss tous les 500 ans) ; avec cette terminologie, Archimède serait de classe 5 000.

Les travaux d'Archimède sont un modèle de clarté, mais la lecture en est difficile: obscurités dans la traduction, différences essentielles dans les modes de pensée. Lorsqu'on a réussi à comprendre, on hésite entre l'admiration et l'exaspération : comment diable fait-il pour savoir que l'ensemble des rationnels n'est pas complet, alors que la notion d'espace métrique complet n'a été introduite par Cauchy qu'en 1820 ? »